\documentstyle[12pt]{article}

\title { Twisted product of Lie groups }
\author { Michael A. Rudkovski }

\date { January 25, 1996 }

\newenvironment{prf}{}{~\hspace{\fill}~$\Box$}
\newenvironment{proof}{\begin{prf}{\bf Proof.}}{\end{prf}}
\newenvironment{example}[1]{\noindent{\bf Example #1.}}{}

\newtheorem{teo}[subsection]{Theorem}
\newtheorem{dest}[subsection]{Proposition}

\newtheorem{remark}[subsection]{Remark}
\newtheorem{defi}[subsection]{Definition}

\begin{document}

\maketitle

\section{Introduction}

\hspace*{\parindent}The constructions of direct and semidirect
products of groups are very frequently used in mathematics.
Firstly, we can construct new groups that inherit certain properties
of initial groups.
Secondly, we can reduce certain "difficult" groups
to several "simple" groups.

One need to note that the semidirect product was appeared as
the generalization of a direct product.

In this article we define the twisted product as
the generalization of the semidirect product of groups.
A twisted product of groups is not always a group.
We will find the necessary and sufficient condition in order that
the twisted product of groups to be a group.
In particular, for two copies of the same group,
the twisted product of group by itself through the action of
inner automorphisms is a group if and only if
the initial group is a metabelian group.

Further we will construct Lie algebra for Lie group of
a twisted product of Lie groups.

After this, we calculate structural constants of
Lie algebra for twisted product of Lie groups
through known structural constants of Lie algebras
of initial Lie groups.

In the case of twisted product of Lie group by itself
by means of the action of inner automorphisms we find
the dependence of the scalar curvature for resulting Lie group
on the scalar curvature for initial Lie group.

\section{Construction of a twisted product of groups}

\hspace*{\parindent}For definition of new operation
we remind some known algebraic concepts.

{\it Action} $\lambda$ of the group $H$ on the group $G$ is the
homomorphism $\lambda : H \rightarrow Aut(G)$.

In particular, for all $g_i \in G$, $h_i \in H$, we have
$\lambda(h_1 h_2)(g) = (\lambda(h_1) \circ \lambda(h_2))(g)$
$= \lambda(h_1)(\lambda(h_2)(g))$,
$\lambda(h)(g_1 g_2) = \lambda(h)(g_1) \lambda(h)(g_2)$,
$\lambda(e_H)(g) = g$,
$\lambda(h)(e_G) = e_G$,
where $e_G$, $e_H$ are identity elements of groups $G$, $H$ accordingly.

For Lie groups $G$ and $H$, the action $\lambda$,
as a map $\lambda: H \times G \rightarrow G$,
is supposed to be smooth; or in another words
$\lambda: H \rightarrow Aut (G)$ is {\it homomorphism of Lie groups}.

{\it Direct product} $G \times H$ of groups $G$ and $H$ is
the set $G \times H$ with the group operation
$(g_1,h_1)(g_2,h_2) = (g_1 g_2,h_1 h_2)$, $g_i \in G$, $h_i \in H$.

{\it Semidirect product} $G \dashv H \equiv G~{\dashv}_{\lambda}~H$
of groups $G$ and $H$ by the action $\lambda : H \rightarrow Aut(G)$
is the set $G \times H$ with the group operation
$(g_1,h_1)(g_2,h_2) = (g_1 \lambda(h_1)(g_2),h_1 h_2)$,
$g_i \in G$, $h_i \in H$.

Here we introduce new concept generalizing the concepts
of direct and semidirect products of groups.
\begin{defi}
Twisted product $G \ast H \equiv G \ast_{\lambda,\mu} H$
of groups $G$ and $H$ by the actions $\lambda : H \rightarrow Aut(G)$
and $\mu : G \rightarrow Aut(H)$ is the set $G \times H$
with the group operation
\begin{equation} \label{mult}
(g_1,h_1)(g_2,h_2) =
(g_1 \lambda(h_1)(g_2), h_1 \mu(g_1)(h_2)), g_i \in G, h_i \in H.
\end{equation}
\end{defi}

In the case of Lie groups we suppose that $\lambda$ and $\mu$ are
homomorphisms of Lie groups.

We notice, that unlike to direct and semidirect products,
a twisted product of groups is not always a group.

\begin{dest} \label{mainth}
A twisted product of groups $G \ast H$ is a group if and only if
the actions $\lambda$ and $\mu$ are satisfied to the condition:
\begin{equation}\label{condition}
\left\{
\begin{array}{r}
\mu ( g) ( h) h^{-1} \in \ker ( \lambda) \\
\lambda ( h) ( g) g^{-1} \in \ker ( \mu)
\end{array} \right.
\end{equation}
for any elements $g \in G$, $h \in H$.
\end{dest}
\begin{proof}
The proof consists of the checking of the group properties
for twisted product of groups.

1) {\sf Existence of identity element.} Obviously, the identity element
of group $G\ast H$ is the element $E := (e_G, e_H)$,
where $e_G$, $e_H$ are identity elements of groups $G$ and $H$ accordingly.

2) { \sf Associativity. } Let us take arbitrary elements
$g_i \in G$, $h_i \in H$.
We shall consider separately left and right parts
of the associativity property:
$$
\biggl( (g_1,h_1)(g_2,h_2) \biggr) (g_3,h_3) =
$$
$$
\biggl( g_1 \lambda (h_1)(g_2) \lambda \Bigl( h_1 \mu(g_1)(h_2) \Bigr) (g_3),
h_1 \lambda (g_1)(h_2) \lambda \Bigl( g_1 \mu(h_1)(g_2) \Bigr) (h_3) \biggr);
$$
$$
(g_1,h_1) \biggl( (g_2,h_2)(g_3,h_3) \biggr) =
$$
$$
\biggl( g_1 \lambda (h_1) \Bigl( g_2 \lambda (h_2) (g_3) \Bigr),
h_1 \mu (g_1) \Bigl( h_2 \mu (g_2) (h_3) \Bigr) \biggr);
$$

By the equating the first coordinates, we have:
$$
g_1 \lambda (h_1)(g_2) \lambda \Bigl( h_1 \mu(g_1)(h_2) \Bigr) (g_3) =
g_1 \lambda (h_1) \Bigl( g_2 \lambda (h_2) (g_3) \Bigr)
\Leftrightarrow
$$
$$
\lambda \Bigl( h_1 \mu(g_1)(h_2) \Bigr) (g_3) =
\lambda (h_1) \Bigl( \lambda (h_2) (g_3) \Bigr)
\Leftrightarrow
$$
$$
[\lambda (h_1) \circ \lambda \Bigl( \mu(g_1)(h_2) \Bigr)] (g_3) =
[\lambda (h_1) \circ \lambda (h_2)] (g_3)
\Leftrightarrow
$$
$$
\lambda \Bigl( \mu(g_1)(h_2) \Bigr)(g_3) =
\lambda (h_2)(g_3)
\Leftrightarrow
$$
$$
\mu (g)(h)h^{-1} \in \ker (\lambda), \mbox{ for any } g \in G, h \in H.
$$

Similarly, by the equating the second coordinates, we will receive
$$
\lambda(h)(g)g^{-1} \in \ker(\mu), \mbox{ for any } g \in G, h \in H.
$$

Thus, for the associativity of twisted product of groups $G \ast H$, it is
necessary and sufficient the condition (\ref{condition}).

For the completion of the proof we shall prove that the condition
(\ref{condition}) is enough for the existence of inverse element.

3) {\sf Existence of inverse element.}
We shall find the inverse element $(g_2,h_2)$
by the considering the equality $(g_1,h_1)(g_2,h_2) = E$.

We have
$(g_1,h_1)(g_2,h_2) =$
$(g_1 \lambda(h_1)(g_2),h_1 \mu (g_1)(h_2)) = (e_G,e_H)$.

Further we shall consider the first coordinate, namely
$$
g_1~\lambda(h_1)(g_2)~=~e_G \Leftrightarrow
\lambda (h_1)(g_2)~=~g_1^{-1} \Leftrightarrow
$$
$$
\Leftrightarrow g_2~=~[\lambda(h_1)]^{-1}(g_1^{-1}) =
\lambda(h_1^{-1})( g_1^{-1}).
$$

Similarly, for the second coordinate we receive
$h_2 = \mu (g_1^{-1})(h_1^{-1})$.

Thus, we have
\begin{equation} \label{back}
(g_1,h_1)^{-1} =
\Bigl(\lambda (h_1^{-1})(g_1^{-1}), \mu (g_1^{-1})(h_1^{-1})\Bigr).
\end{equation}

Then the condition that right-inverse element must be also
left-inverse element give us
$$
\Bigl(\lambda (h_1^{-1})(g_1^{-1}),\mu (g_1^{-1})(h_1^{-1})\Bigr)(g_1,h_1) =
$$
$$
\biggl(\lambda (h_1^{-1})(g_1^{-1})
\lambda \Bigl( \mu (g_1^{-1})(h_1^{-1})\Bigr)(g_1),
\mu (g_1^{-1})(h_1^{-1})
\mu \Bigl(\lambda(h_1^{-1})(g_1^{-1})\Bigr)(h_1)\biggr) = E
$$

By the equating the first coordinates, we receive
$$
\lambda (h_1^{-1})(g_1^{-1})
\lambda \Bigl( \mu (g_1^{-1})(h_1^{-1})\Bigr)(g_1) = e_G \Leftrightarrow
$$
$$
\lambda \Bigl( \mu (g_1^{-1})(h_1^{-1})\Bigr)(g_1) =
[\lambda (h_1^{-1})(g_1^{-1})]^{-1} = \lambda (h_1^{-1})(g_1)
$$

For given condition, as easily to see, it is enough
to have the first part of the condition (\ref{condition}), namely
$$
\mu (g)(h)h^{-1} \in \ker (\lambda), \mbox{ for any } g \in G, h \in H.
$$

Similarly we receive that for the equality of the second coordinates
it is enough to have the second part
of the condition (\ref{condition}), namely
$$
\lambda (h)(g)g^{-1} \in \ker (\mu), \mbox{ for any } g \in G, h \in H.
$$

Thus, the existence of the inverse element
in the twisted product of groups $G \ast H$
follows from the condition (\ref{condition}).
\end{proof}

\begin{remark} \label{rem_smooth}
The smoothness of homomorphisms $\lambda$ and $\mu$
implies the smoothness of the multiplication (\ref{mult}) relative to
the initial smooth structure of direct product of smooth manifolds,
even if the twisted product is not a Lie group.
\end{remark}

\begin{dest} \label{lie}
If twisted product of Lie groups is a group,
then the group is a Lie group.
\end{dest}
\begin{proof}
The smoothness of the group multiplication follows
from the previous remark.
Smoothness of the operation of the taking of inverse element follows
directly from the formula (\ref{back}).
\end{proof}

\begin{dest} \label{mn}
If $G$ and $H$ are two copies of the same group $M \equiv G \equiv H$,
and $\lambda$, $\mu$ are the actions of group $M$ on itself by inner
automorphisms, i.e.
$\lambda (h)(g) = hgh^{-1}$, $\mu (g)(h) = ghg^{-1}$, $g$, $h \in M$,
then the twisted product of groups $G \ast H \equiv M \ast M$
by the actions $\lambda$, $\mu$ is a group if and only if
the initial group $M \equiv G \equiv H$ is metabelian
(2-step nilpotent) group.
\end{dest}
\begin{proof}
We consider the condition (\ref{condition}):
$$
\left\{
\begin{array}{r}
\mu ( g) ( h) h^{-1} \in \ker ( \lambda) \\
\lambda ( h) ( g) g^{-1} \in \ker ( \mu)
\end{array} \right.
\Leftrightarrow
$$
$$
\left\{
\begin{array}{r}
ghg^{-1}h^{-1} \in \ker ( \lambda) \\
hgh^{-1}g^{-1} \in \ker ( \mu)
\end{array} \right.
$$

As far as $g \in M$, $h \in M$ are any elements,
then for arbitrary $x \in M$,
we have the chain of the equivalences:
$$
\Bigr( (ghg^{-1}h^{-1})x(hgh^{-1}g^{-1}) \Bigl)x^{-1} = e_M
\Leftrightarrow
$$
$$
[g,h]x[g,h]^{-1}x^{-1} = e_M \Leftrightarrow
[[g,h],x] = e_M, g,h,x \in M.
$$
\end{proof}

\begin{remark} \label{remar}
If in twisted product of group on itself (see proposition \ref{mn})
we take outer automorphisms,
then the requirement of metabelian condition is not obligatory
for group condition.
\end{remark}

\begin{example}{1}
As an example we shall consider the twisted product
of the 3-dimensional group $E(2)$ proper movements
of 2-dimensional Euclidean space on itself by the actions
$\lambda (A,\xi)(B,\eta) = \mu (A,\xi)(B,\eta) = (B, A \eta)$,
where $A$, $B \in SO(2)$, $\xi$, $\eta \in P$,
and $P \cong R^2$ is the group of parallel translations on Euclidean plane.
(Known decomposition $E (2)$ into the semidirect product $SO (2) \vdash P$
is also considered. The multiplication in $E (2)$ is defined as follows:
$(A,\xi)(B,\eta) = (AB, \xi + A \eta)$,
$A$,$B \in SO(2)$, $\xi$, $\eta \in P$).
One easily checks up that these actions $\lambda$, $\mu$ satisfy
to the condition (\ref{condition}), so $E(2) \ast E(2)$ is a group, though
the solvable group $E(2)$ is not metabelian group.
\end{example}

\section{Lie algebra for twisted product of Lie groups}

\hspace*{\parindent}For Lie groups $G$ and $H$
by $L(G)$, $L(H)$ we shall denote its Lie algebras.

We shall assume that the twisted product of Lie groups $G \ast H$
by actions $\lambda$, $\mu$ satisfies to the condition (\ref{condition}).
Then according to the propositions \ref{mainth} and \ref{lie}
$G \ast H$ are Lie groups.

Obviously, as a vector space, Lie algebra $L(G \ast H)$ of Lie group of
twisted product of Lie groups $G \ast H$
is represented as direct sum of vector spaces $L(G) \bigoplus L(H)$.
Next we calculate Lie bracket in Lie algebra
$L(G \ast H)$.

Let $\lambda~:~H~\rightarrow~Aut(G)$, $\mu~:~G~\rightarrow~Aut(H)$ are
appropriate actions of Lie groups. Then its differentials
$d {\lambda (h)}_{e} \in Aut(L(G))$ at the identity element
of Lie group $G$ and
$d {\mu (g)}_{e} \in Aut(L(H))$ at the identity element
of Lie group $H$ are automorphisms of Lie algebras
$L(G)$, $L(H)$ accordingly.

Thus, $d {\lambda ( \cdot )}_{e} : H \rightarrow Aut L(G)$,
$d {\mu ( \cdot )}_{e} : G \rightarrow Aut L(H)$
are homomorphisms of Lie groups into linear groups.
Its differentials at the identity element of groups $H$, $G$
accordingly we shall denote:
$L := d {\lambda ( {\cdot}_{e} )}_{e} : L(H) \rightarrow Der L(G)$,
$M := d {\mu ( {\cdot}_{e} )}_{e} : L(G) \rightarrow Der L(H)$,
where $Der L(G)$, $Der L(H)$ are Lie algebras of the derivations
for Lie algebras $L(G)$, $L(H)$.

Let
$g_i = \exp tX_i,~ h_i = \exp tY_i;~ X_i \in L(G),~ Y_i \in L(H);~ i = 1,2$,
then $(g_1,h_1)(g_2,h_2) = $
$\Bigl( \exp tX_1~ \lambda ( \exp tY_1)( \exp tX_2),$
$\exp tY_1~ \mu ( \exp tX_1)( \exp tY_2) \Bigr) = $
$\Bigl( \exp tX_1~ \exp ( {e}^{L(tY_1)}(tX_2)), $
$\exp tY_1~ \exp ( {e}^{M(tX_1)}(tY_2)) \Bigr) $,
where ${e}^A = {\displaystyle \sum_{k=0}^{\infty}} \frac{A^k}{k!}$
is the exponential map for the group of automorphisms of the
appropriate vector spaces ($L(G)$ or $L(H)$).

In fact we consider only $\lambda~:~H~\rightarrow~Aut(G)$.
We have the following commutative diagrams:

$$
\begin{array}{rcl}
tX_2 \in L(G) &
\stackrel{d {\lambda ( \exp tY_1)}_e}{\longrightarrow} & L(G) \\
{\scriptstyle \exp} \downarrow & {} & \downarrow {\scriptstyle \exp} \\
G & \stackrel{\lambda ( \exp tY_1)}{\longrightarrow} & G
\end{array}
$$

$$
\begin{array}{rcl}
tY_1 \in L(H) &
\stackrel{d {\lambda ( {\cdot}_e )}_e \equiv L}{\longrightarrow} &
Der L(G) \\
{\scriptstyle \exp} \downarrow & {} &
\downarrow {\scriptstyle \exp \equiv {e}} \\
H & \stackrel{d {\lambda ( \cdot )}_e}{\longrightarrow} & Aut L(G)
\end{array}
$$

From this we receive:
$$
\lambda ( \exp tY_1)( \exp tX_2)
= \exp ( d {\lambda ( \exp tY_1)}_e (tX_2) )
= \exp ( {e}^{L(tY_1)}(tX_2) ).
$$

As a result we have:
$$
(g_1,h_1)(g_2,h_2) =
\biggl( \exp \Bigl( tX_1 + tX_2 + \frac{t^2}{2}[X_1,X_2] +
t^2 L(Y_1)(X_2) + O(t^3) \Bigr),
$$
$$
\exp \Bigl( tY_1 + tY_2 + \frac{t^2}{2}[Y_1,Y_2] +
t^2 M(X_1)(Y_2) + O(t^3) \Bigr) \biggr)
=: \exp \xi(t), \xi(t) \in L(G \ast H);
$$

$$
(g_2,h_2)(g_1,h_1) =
\biggl( \exp \Bigl( tX_2 + tX_1 + \frac{t^2}{2}[X_2,X_1] +
t^2 L(Y_2)(X_1) + O(t^3) \Bigr),
$$
$$
\exp \Bigl( tY_2 + tY_1 + \frac{t^2}{2}[Y_2,Y_1] +
t^2 M(X_2)(Y_1) + O(t^3) \Bigr) \biggr)
=: \exp \eta(t), \eta(t) \in L(G \ast H).
$$

Further we shall consider the commutator of elements
$(g_1,h_1)$ and $(g_2,h_2)$ in the group $G \ast H$.
We have
$$
[(g_1,h_1),(g_2,h_2)] =
(g_1,h_1)(g_2,h_2) {\Bigl( (g_2,h_2)(g_1,h_1) \Bigr)}^{-1} =
$$
$$
\exp \xi(t) {\Bigl( \exp \eta(t) \Bigr)}^{-1} =
\exp \xi(t) \exp \Bigl( - \eta(t) \Bigr) =
\exp \Bigl( \xi(t) - \eta(t) + O(t^3) \Bigr).
$$

Then for $Z_i = (X_i,Y_i) \in L(G \ast H)$, $i = 1,2$ we receive
the formula for the calculations of Lie bracket in Lie algebra
$L(G \ast H)$ of Lie group of twisted product of Lie groups:
$$
[Z_1,Z_2] =
{\frac{1}{2} \frac{d^2}{dt^2}|}_{t=0}
\Bigl( \xi(t) - \eta(t) + O(t^3) \Bigr) =
$$
\begin{equation} \label{mainform}
\Bigl( [X_1,X_2] + L(Y_1)(X_2) - L(Y_2)(X_1),
[Y_1,Y_2] + M(X_1)(Y_2) - M(X_2)(Y_1) \Bigr).
\end{equation}

\begin{example}{2}
Now we illustrate the technique to the application of these results.
We find well known Lie algebra. It is a Lie algebra of
Heisenberg group ${\Gamma}_M$,
i.e. the group of upper triangular $3 \times 3$ matrixes with
identity at a main diagonal, under the matrix product.

Associate with each element
$$
\left( \begin{array}{ccc}
1 & c & b \\
0 & 1 & a \\
0 & 0 & 1
\end{array} \right) \in {\Gamma}_M
$$
the element $(a,b,c) \in R^3$ and introduce on $R^3$
the product operation
$$
(a_1,b_1,c_1) \cdot (a_2,b_2,c_2) = (a_1+a_2,b_1+b_2+c_1 a_2,c_1+c_2).
$$

Then we shall receive the group $\Gamma = ( R^3, \cdot)$,
which isomorphic to the group ${\Gamma}_M$.
Let us represent the group $\Gamma$ in the form
of semidirect product $R^2 \dashv R^1$.
For this purpose we shall define the action $\lambda$
of the group $(R^1,+)$ on the group $(R^2,+)$ as following:
$$
\lambda : R^1 \ni h \longrightarrow \phi_h :=
(R^2 \ni (g^1,g^2) \rightarrow (g^1, g^2 + h g^1) \in R^2) \in Aut(R^2)
$$

As far as a semidirect product is a partial case of a twisted product,
then to the calculation of a Lie bracket it is possible to apply
expression (\ref{mainform}) received earlier.

Let $X=(x^1,x^2)$, $Y=(y)$ are generating vectors of one-parameter
subgroups $g(t), h(t)$ accordingly.

Then we have
$$
\begin{array}{c}
L(Y)(X) = {\frac{d^2}{dt ds}}|_{t=0,s=0}
\lambda \biggl( h(t)\biggr)\biggl( g(s)\biggr) = \\
{\frac{d^2}{dt ds}}|_{t=0,s=0}
\biggl( g^{1}(s), g^{2}(s)+h(t)g^{1}(s) \biggr) = \\
\biggl( 0, {\frac{d}{dt}}|_{t=0}h(t) ~ {\frac{d}{dt}}|_{t=0}g^{1}(s) \biggr)
= (0,y x^1)
\end{array}
$$

Applying the map $L$, mentioned earlier, we receive
for vectors
$Z_i=(X_i,Y_i)$, $X_i=(x^1_i,x^2_i)$, $Y_i=(y_i)$,
$i=1,2$,
the equation:
\begin{equation} \label{primer}
[Z_1,Z_2] = [(x^1_1,x^2_1,y_1),(x^1_2,x^2_2,y_2)] =
(0,y_1 x^1_2 - y_2 x^1_1,0)
\end{equation}

Thus, the Lie algebra of Lie group $\Gamma$ is
the vector space $R^3$
with Lie bracket calculated by the formula (\ref{primer}).
\end{example}

\begin{dest} \label{pm}
If $G$ and $H$ are two copies of the same metabelian
Lie group $M \equiv G \equiv H$, and $\lambda$,  $\mu$ are
actions of the Lie group $M$ on itself by inner automorphisms, i.e.
$\lambda (h)(g) = hgh^{-1}$, $\mu (g)(h) = ghg^{-1}$, $g,h \in M$,
then twisted product of Lie groups
$G~\ast_{\lambda,\mu}~H \equiv M~\ast_{\lambda,\mu}~M$
is also metabelian Lie group.
\end{dest}
\begin{proof}
We take arbitrary vectors $Z_1, Z_2, Z_3 \in L(M \ast M)$.
As far as $\lambda$, $\mu$ are the actions by inner automorphisms,
then for arbitrary vectors $X,Y \in L(M)$,
$L(X)(Y) \equiv M(X)(Y) = ad_X Y = [X,Y]$.
Therefore the Lie bracket (\ref{mainform}) accepts a form:
\begin{equation} \label{pmainform}
[Z_1,Z_2] =
\Bigl( [X_1,X_2] + [Y_1,X_2] - [Y_2,X_1],
[Y_1,Y_2] + [X_1,Y_2] - [X_2,Y_1] \Bigr).
\end{equation}

It is known that a Lie group $M$ is metabelian if and only if
its Lie algebra $L(M)$ is metabelian.
It means by the definition that for any vectors $X_1,X_2,X_3~\in~L(M)$,
we have the equality $[[X_1,X_2],X_3] = 0$.

Hence it is easily to check up that $[[Z_1,Z_2],Z_3] = 0$,
i.e. the Lie algebra $L(M \ast M)$ is metabelian.
Therefore, Lie group $M \ast M$ is metabelian.
\end{proof}

\section{Structural constants of Lie algebra for Lie Group of
twisted product of Lie groups}

\hspace*{\parindent}Let $G$ and $H$ are Lie groups.
By $L(G)$, $L(H)$ we shall denote its Lie algebras.

We shall assume that a twisted product of Lie groups $G \ast H$
by the actions $\lambda$, $\mu$ is satisfied
to the condition (\ref{condition}).
Then according to the propositions \ref{mainth} and \ref{lie}
$G \ast H$ is a Lie group.

As was calculated above, the Lie algebra $L(G \ast H)$ of Lie group
$G \ast H$ is a vector space $L(G) \bigoplus L(H)$ with Lie bracket
defined by the formula (\ref{mainform}).

Let $n = \dim L(G)$, $m = \dim L(H)$, ${ \{ e_i^G \} }_{i=1}^n$
is a basis of Lie algebra $L(G)$, and
${ \{ e_i^H \} }_{i=1}^m$ is a basis of Lie algebra $L(H)$.
Let ${0}^G$, ${0}^H$ are zero elements of Lie algebras
$L(G)$ and $L(H)$ accordingly.
We shall introduce the basis of the Lie algebra $L(G \ast H)$ as follows:
\begin{equation} \label{base}
E_i = \left\{
\begin{array}{ll}
(e_i^G,{0}^H), & \mbox{ for }i \leq n \\
({0}^G,e_{i-n}^H), & \mbox{ for }i > n
\end{array} \right.
\end{equation}

Let us denote the projection of the vector $X \in L(G \ast H)$ onto
the first factor $L(G)$ as $X^1 = {\Pr}_{L(G)}X$, and onto the second
factor $L(H)$ as $X^2 = {\Pr}_{L(H)}X$.

We need to calculate structural constants
of Lie algebra $L(G \ast H)$,
i.e. the numbers $c_{ij}^k$,
such that $[E_i,E_j] = {\displaystyle \sum_{k=1}^{n+m}} c_{ij}^k E_k$,
$1 \leq i \leq n+m$, $1 \leq j \leq n+m$.
We shall put that structural constants of Lie algebras
$L(G)$ and $L(H)$ in the bases
${ \{ e_i^G \} }_{i=1}^n$ and ${ \{ e_i^H \} }_{i=1}^m$
are equal accordingly ${ \{ {\beta}_{ij}^k \} }_{i,j,k=1}^n$ and
${ \{ {\gamma}_{ij}^k \} }_{i,j,k=1}^m$.

\begin{dest}
In the notations above, the structural constants $c_{ij}^k$
of Lie algebra $L(G \ast H)$ of twisted product of Lie groups
$G \ast H$ equal to
\begin{equation}\label{star}
c^k_{ij} = \left\{
\begin{array}{rl}
\mbox{ for }i \leq n, j \leq n & \left\{ \begin{array}{lr} \beta^k_{ij}, &
\mbox{ for }k \leq n \\[2mm] 0, &
\mbox{ for }k > n \end{array} \right. \\ \vspace{4pt}
\mbox{ for }i > n, j > n & \left\{ \begin{array}{lr} 0, &
\mbox{ for }k \leq n \\[2mm] \gamma^{k-n}_{(i-n) (j-n)}, &
\mbox{ for }k > n \end{array} \right. \\ \vspace{4pt}
\mbox{ for }i \leq n, j > n &
\left\{ \begin{array}{lr} -< L(e^H_{j-n})(e^G_i),e^G_k>, &
\mbox{ for }k \leq n \\[2mm] < M(e^G_i)(e^H_{j-n}),e^H_{k-n}>, &
\mbox{ for }k > n \end{array} \right. \\ \vspace{4pt}
\mbox{ for }i > n, j \leq n &
\left\{ \begin{array}{lr} < L(e^H_{i-n})(e^G_j),e^G_k>, &
\mbox{ for }k \leq n \\[2mm] -< M(e^G_j)(e^H_{i-n}),e^H_{k-n}>, &
\mbox{ for }k > n \end{array} \right. \\
\end{array} \right.
\end{equation}
\end{dest}
\begin{proof}
From the formula (\ref{mainform}) we have
$$
\begin{array}{c}
[E_i,E_j] = \bigl( L(E^2_i)(E^1_j) -
L(E^2_j)(E^1_i) + [E_i^1,E_j^1], M(E^1_i)(E^2_j) -\\
M(E^1_j)(E^2_i) + [E_i^2,E_j^2] \bigr).
\end{array}
$$
Now we consider four cases, namely

1) $i \leq n$, $j \leq n$. With regard
to the basis ${ \{ E_i \} }_{i=1}^{n+m}$ in (\ref{base}) we receive
$$
[E_i,E_j] = ([e_i^G,e_j^G],{0}^H).
$$
Thus, we have
$$
c^k_{ij} = \left\{
\begin{array}{ll}
{\beta}^k_{ij}, & \mbox{ for }k \leq n, \\
0, & \mbox{ for }k > n.
\end{array} \right.
$$

2) $i > n$, $j > n$. Similarly we receive
$$
c^k_{ij} = \left\{
\begin{array}{ll}
0, & \mbox{ for }k \leq n,\\
{\gamma}^{k-n}_{(i-n) (j-n)}, & \mbox{ for }k > n.
\end{array} \right.
$$

3) $i \leq n$, $j > n$. In this case we have
$$
\begin{array}{c}
[E_i,E_j] = -< L(e^H_{j-n})(e^G_i),M(e^G_i)(e^H_{j-n})> =
{ \displaystyle \sum_{k=1}^{n+m}} c_{ij}^k E_k = \\
{ \displaystyle \sum_{k=1}^{n}} c_{ij}^k (e_k^G,{0}^H) +
{ \displaystyle \sum_{k=1}^{m}} c_{ij}^{k+n} ({0}^H,e_k^H)
\end{array}
$$
We have
$$
- L(e^H_{j-n})(e^G_i) = { \displaystyle \sum_{k=1}^{n}} c_{ij}^k e_k^G,
$$
$$
M(e^G_i)(e^H_{j-n}) = { \displaystyle \sum_{k=1}^{m}} c_{ij}^{k+n} e_k^H.
$$
We receive in this case
$$
c^k_{ij} = \left\{
\begin{array}{rl}
- <L(e^H_{j-n})(e^G_i), e_k^G>, & \mbox{ for }k \leq n,\\
<M(e^G_i)(e^H_{j-n}), e_{k-n}^H>, & \mbox{ for }k > n.
\end{array} \right.
$$

4) $i > n$, $j \leq n$.
$$
[E_i,E_j] = < L(e^H_{i-n})(e^G_j), - M(e^G_j)(e^H_{i-n})> =
{ \displaystyle \sum_{k=1}^{n}} c_{ij}^k (e_k^G,{0}^H) +
{ \displaystyle \sum_{k=1}^{m}} c_{ij}^{k+n} ({0}^H,e_k^H)
$$
We have
$$
L(e^H_{i-n})(e^G_j) = { \displaystyle \sum_{k=1}^{n}} c_{ij}^k e_k^G,
$$
$$
- M(e^G_j)(e^H_{i-n}) = { \displaystyle \sum_{k=1}^{m}} c_{ij}^{k+n} e_k^H.
$$
We receive in this case
$$
c^k_{ij} = \left\{
\begin{array}{rl}
<L(e^H_{i-n})(e^G_j), e_k^G>, & \mbox{ for }k \leq n,\\
- <M(e^G_j)(e^H_{i-n}), e_{k-n}^H>, & \mbox{ for }k > n.
\end{array} \right.
$$
\end{proof}

\begin{dest}
Structural constants ${ \{ c_{i,j}^k \} }_{i,j,k=1}^{2n}$
of Lie algebra for twisted product of
metabelian Lie group $M$ by itself $M \ast M$ by means of the actions of
inner automorphisms (see proposition \ref{mn}) are expressed in
structural constants of the initial Lie algebra $L(M)$ as follows:
\begin{equation} \label{property}
c^k_{ij} = \left\{
\begin{array}{rl}
\mbox{ for }i \leq n, j \leq n &
\left\{ \begin{array}{lr} \alpha^k_{ij}, &
\mbox{ for }k \leq n \\ 0, &
\mbox{ for }k > n \end{array} \right. \\ \vspace{4pt}
\mbox{ for }i > n, j > n &
\left\{ \begin{array}{lr} 0, &
\mbox{ for }k \leq n \\ \alpha^k_{ij}, &
\mbox{ for }k > n \end{array} \right. \\ \vspace{4pt}
\mbox{ for }i \leq n, j > n &
\left\{ \begin{array}{lr} \alpha^k_{i (j-n)}, &
\mbox{ for }k \leq n \\[2mm] \alpha^{k-n}_{i (j-n)}, &
\mbox{ for }k > n \end{array} \right. \\ \vspace{4pt}
\mbox{ for }i > n, j \leq n &
\left\{ \begin{array}{lr} \alpha^k_{(i-n) j}, &
\mbox{ for }k \leq n \\[2mm] \alpha^{k-n}_{(i-n) j}, &
\mbox{ for }k > n \end{array} \right.
\end{array} \right.
\end{equation}
\end{dest}
\begin{proof}
As far as $\lambda$, $\mu$ are the actions of
inner automorphisms of Lie group $M$ by itself,
then for arbitrary vectors $X,Y \in L(M)$,
$L(X)(Y) \equiv M(X)(Y) = ad_X Y = [X,Y]$.

1),2) The proof for $i \leq n$, $j \leq n$ and
for $i > n$, $j > n$ follows directly from the formula (\ref{star}).

3) For $i \leq n$, $j > n$ we have
$$
k \leq n \Rightarrow c_{ij}^k =
-<L(e^H_{j-n})(e^G_i), e_k^G> = -<[e_{j-n},e_i],e_k> =
\alpha^k_{i (j-n)}
$$
$$
k > n  \Rightarrow c_{ij}^k =
<M(e^G_i)(e^H_{j-n}), e_{k-n}^H> = <[e_i,e_{j-n}],e_{k-n}> =
\alpha^{k-n}_{i (j-n)}
$$

4) For $i > n$, $j \leq n$ we have
$$
k \leq n \Rightarrow c_{ij}^k =
<L(e^H_{i-n})(e^G_j), e_k^G> = <[e_{i-n},e_j],e_k> =
\alpha^k_{(i-n) j}
$$
$$
k > n  \Rightarrow c_{ij}^k =
- <M(e^G_j)(e^H_{i-n}), e_{k-n}^H> = - <[e_j,e_{i-n}],e_{i-n}],e_k> =
\alpha^{k-n}_{(i-n) j}
$$
\end{proof}

\section{Scalar curvature of Lie group of
the twisted product of Lie group by itself}

\hspace*{\parindent}Let $G \equiv H \equiv M$ are two copies
of the same Lie group $M$.
We consider their twisted product $G \ast H \equiv M \ast M$, received
by the actions of inner automorphisms
$\lambda$, $\mu : M \ni h  \longrightarrow h g h^{-1} \in Aut(M)$.

We assume further in all cases that the given actions
$\lambda$, $\mu$ satisfy to the condition (\ref{condition}),
hence, according to the propositions \ref{mainth}, \ref{lie},
twisted product of Lie group by itself
$M \ast_{\lambda, \mu} M$ is a Lie group.

\begin{teo} \label{scalcurv}
Let $\rho'$ is the scalar curvature of Lie group
$M \ast_{\lambda, \mu} M$,
where $\lambda$, $\mu$ are the actions of inner automorphisms,
and $\rho$ is the scalar curvature of Lie group $M$.
Then we have the equation:
\begin{equation} \label{scalequat}
\rho' = 6 \rho
\end{equation}
\end{teo}
\begin{proof}
As far as the twisted product $M \ast M$ by the actions of inner
automorphisms is a group, then it follows from the proposition \ref{mn}
that group $M$ is a metabelian group.

Since $M$ is metabelian, then it follows from
the proposition \ref{pm} that Lie group $M \ast M$
is also metabelian.

Because Lie group $M$ is metabelian,
then scalar curvature $\rho$ of Lie group $M$ is
calculated as follows (see \cite{M}, \cite[page 11]{B})
\begin{equation} \label{scal}
\rho = - \frac{1}{4} {\displaystyle \sum_{i=1}^n \sum_{k=1}^n}
{\Vert [e_i,e_k] \Vert}^2
\end{equation}

As far as Lie group $M \ast M$ is metabelian, we have similarly
\begin{equation} \label{pscal}
\rho' =
- \frac{1}{4} {\displaystyle \sum_{i=1}^{2n} \sum_{k=1}^{2n}}
{\Vert [E_i,E_k] \Vert}^2
\end{equation}

Thus, from the expression (\ref{property}) we receive
$$
\begin{array}{c}
\rho' =
- \frac{1}{4} {\displaystyle \sum_{i=1}^{2n} \sum_{k=1}^{2n}}
{\Vert [E_k,E_i] \Vert}^2 =
- \frac{1}{4} {\displaystyle \sum_{i=1}^{n} \sum_{k=1}^{n}}
{\Vert [E_k,E_i] \Vert}^2
- \frac{1}{4} {\displaystyle \sum_{i=1}^{n} \sum_{k=n+1}^{2n}}
{\Vert [E_k,E_i] \Vert}^2 - \\
- \frac{1}{4} {\displaystyle \sum_{i=n+1}^{2n} \sum_{k=1}^{n}}
{\Vert [E_k,E_i] \Vert}^2
- \frac{1}{4} {\displaystyle \sum_{i=n+1}^{2n} \sum_{k=n+1}^{2n}}
{\Vert [E_k,E_i] \Vert}^2 =
- \frac{1}{4} {\displaystyle \sum_{i=1}^{n}
\sum_{k=1}^{n} \sum_{j=1}^{2n}} (c_{ki}^j)^2 - \\
- \frac{1}{4} {\displaystyle \sum_{i=1}^{n}
\sum_{k=n+1}^{2n} \sum_{j=1}^{2n}} (c_{ki}^j)^2
- \frac{1}{4} {\displaystyle \sum_{i=n+1}^{n}
\sum_{k=1}^{n} \sum_{j=1}^{2n}} (c_{ki}^j)^2
- \frac{1}{4} {\displaystyle \sum_{i=n+1}^{n}
\sum_{k=n+1}^{2n} \sum_{j=1}^{2n}} (c_{ki}^j)^2 = \\
- \frac{1}{4} {\displaystyle \sum_{i=n}^{n}
\sum_{k=1}^{n} \sum_{j=1}^{n}} ({\alpha}_{ki}^j)^2
- \frac{1}{2} {\displaystyle \sum_{i=1}^{n}
\sum_{k=1}^{n} \sum_{j=1}^{n}} ({\alpha}_{ki}^j)^2
- \frac{1}{2} {\displaystyle \sum_{i=1}^{n}
\sum_{k=1}^{n} \sum_{j=1}^{n}} ({\alpha}_{ki}^j)^2 - \\
- \frac{1}{4} {\displaystyle \sum_{i=1}^{n}
\sum_{k=1}^{n} \sum_{j=1}^{n}} ({\alpha}_{ki}^j)^2 =
- \frac{3}{2} {\displaystyle \sum_{i=1}^{n}
\sum_{k=1}^{n} \sum_{j=1}^{n}} ({\alpha}_{ki}^j)^2 =
6\rho.
\end{array}
$$
\end{proof}

\begin{example}{3}
As an example we shall take the twisted product ${\Gamma} * {\Gamma}$
of Heizenberg group by itself under the actions of inner automorphisms
$$
\lambda, \mu : \Gamma \ni g \mapsto
(\Gamma \ni h \mapsto ghg^{-1} \in \Gamma) \in Aut(\Gamma).
$$

As far as the Heizenberg group $\Gamma$ is metabelian,
then the twisted product ${\Gamma} * {\Gamma}$ is a group.

In coordinates the actions $\lambda$, $\mu$ have a form:
$$
\lambda (A_1,B_1,C_1) (A_2,B_2,C_2) =
\mu (A_1,B_1,C_1) (A_2,B_2,C_2) =
$$
$$
= (A_2, B_2 + C_1 A_2 - C_2 A_1, C_2).
$$

Thus, the operators $L$, $M$ have a form:
$$
L (a_1,b_1,c_1) (a_2,b_2,c_2) =
M (a_1,b_1,c_1) (a_2,b_2,c_2) =
(0, c_1 a_2 - c_2 a_1, 0).
$$

We shall remind that in the example 1 we calculate the Lie bracket
(\ref{primer}) of Lie algebra of Heisenberg group $\Gamma$.
Namely, the Lie algebra $\Gamma$
of Lie group $\Gamma \cong R^2 \dashv R^1$
is a vector space $R^3$ with Lie bracket defined by the formula:
$$
[Z_1,Z_2] = [(x^1_1,x^2_1,y_1),(x^1_2,x^2_2,y_2)]
= (0,y_1 x^1_2 - y_2 x^1_1,0)
$$
where $Z_i=(X_i,Y_i) \in \Gamma \cong R^2 \dashv R^1$,
$X_i=(x^1_i,x^2_i) \in R^2$, $Y_i=(y_i) \in R^1$, $i = 1,2$.

Further from the formula (\ref{mainform}) we receive
that Lie bracket in $L({\Gamma} * {\Gamma})$ accepts a form:
$$
[Z_1,Z_2] =
\Bigl( [X_1,X_2] + L(Y_1)(X_2) - L(Y_2)(X_1),
$$
$$
[Y_1,Y_2] + M(X_1)(Y_2) - M(X_2)(Y_1) \Bigr) =
$$
$$
(0, a^1_1 c^1_2 - a^1_2 c^1_1 + a^2_1 c^1_2
- a^1_2 c^2_1 - a^2_2 c^1_1 + a^1_1 c^2_2, 0,
$$
$$
0, a^2_1 c^2_2 - a^2_2 c^2_1 + a^1_1 c^2_2
- a^2_2 c^1_1 - a^1_2 c^2_1 + a^2_1 c^1_2, 0)
$$

We shall choose orthonormal basis $\{e_i\}$ in all
copies of $\Gamma$ as follows:
$e_1 = (1,0,0)$, $e_2 = (0,1,0)$, $e_3 = (0,0,1)$.

Since, as a vector space,
$L(\Gamma \ast \Gamma) \equiv L(\Gamma) \oplus L(\Gamma)$,
then in Lie algebra $L(\Gamma \ast \Gamma)$ naturally arises
orthonormal basis
$$
\{E_1,E_2,E_3,E_4,E_5,E_6\} =
\{e_1 \oplus 0, e_2 \oplus 0, e_3 \oplus 0,
0 \oplus e_1, 0 \oplus e_2, 0 \oplus e_3\}
$$

Thus, Lie bracket of Lie algebra
$L(\Gamma \ast \Gamma)$
for the basis vectors
${\{E_i\}}_{i=1}^6$ accepts the values
(here we show only nonzero values):
$$
\begin{array}{lcr}
~[E_1,E_3] & = & - E_2, \\
~[E_1,E_6] & = & - E_2 - E_5, \\
~[E_4,E_3] & = & - E_2 - E_5, \\
~[E_4,E_6] & = & - E_5;
\end{array}
$$

In other words, by the formula (\ref{star}) of structural
constants of Lie algebra for twisted product of Lie groups
we receive the structural constants of Lie algebra
$L(\Gamma \ast \Gamma)$
(here we show only nonzero structural constants):
$$
\alpha_{13}^2=-1, \alpha_{43}^2=-1, \alpha_{43}^5=-1,
\alpha_{16}^2=-1, \alpha_{16}^5=-1, \alpha_{46}^5=-1,
$$
$$
\alpha_{31}^2=1, \alpha_{34}^2=1, \alpha_{34}^5=1,
\alpha_{61}^2=1, \alpha_{61}^5=1, \alpha_{64}^5=1;
$$

We have following formula of sectional curvatures of Lie group
$G$ by the structural constants
of Lie algebra $L(G)$ (see \cite[page 295]{M})
\begin{equation} \label{sectcurv}
k_{ij} =
\sum_{k=1}^{n} ( \frac{1}{2} \alpha_{ij}^k
( - \alpha_{ij}^k + \alpha_{jk}^i + \alpha_{ki}^j)
- \frac{1}{4} ( \alpha_{ij}^k - \alpha_{jk}^i + \alpha_{ki}^j)
( \alpha_{ij}^k + \alpha_{jk}^i - \alpha_{ki}^j)
- \alpha_{ki}^i \alpha_{kj}^j ),
\end{equation}
$$
\mbox{where}~ i,j = 1..n, n = dim L(G).
$$

Then by direct calculation we receive a matrix of sectional curvatures
(the sectional curvature in 2-dimensional direction
stretched on the basis vectors $E_i$, $E_j$
stands on i-th column and j-th row)
of the Lie group $\Gamma \ast \Gamma$:
$$
\left [\begin {array}{cccccc} 0&1/2&-3/4&0&1/4&-3/2
\\\noalign{\medskip}1/2&0&1/2&1/4&0&1/4\\\noalign{\medskip}-3/4&1/2
&0&-3/2&1/4&0\\\noalign{\medskip}0&1/4&-3/2&0&1/2&-3/4
\\\noalign{\medskip}1/4&0&1/4&1/2&0&1/2\\\noalign{\medskip}-3/2&1/4
&0&-3/4&1/2&0\end {array}\right ]
$$

From this we can easily find that the scalar curvature of the group
$\Gamma \ast \Gamma$ is equal $-3$.

Next we shall calculate the scalar curvature
of initial Lie group $\Gamma$.

From the formula (\ref{primer}) we can easily receive
that for basis vectors ${\{e_i\}}_{i=1}^3$
the Lie bracket of Lie algebra $L(\Gamma)$ accepts
the single nonzero value:
$$
[e_1,e_3] = - e_2.
$$

Then we receive a matrix of sectional curvatures
of the Lie group $\Gamma \ast \Gamma$:
$$
\left [\begin {array}{ccc} 0&1/4&-3/4\\\noalign{\medskip}1/4&0&1
/4\\\noalign{\medskip}-3/4&1/4&0\end {array}\right ]
$$

Thus, the scalar curvature
of Lie group $\Gamma$ is equal $- \frac{1}{2}$.

We see that the scalar curvature of initial Lie group $\Gamma$
and the scalar curvature of $\Gamma \ast \Gamma$
satisfy to equation (\ref{scalequat}).
\end{example}

\begin{example}{4}
As a next example we shall take, considered above in the example 1,
the twisted product of group $E (2) \cong SO (2) \vdash P$ of proper
movements of the 2-dimensional Euclidean space by itself under the actions
$\lambda (A,\xi)(B,\eta) = \mu (A,\xi)(B,\eta) = (B, A \eta)$.
Here $A$, $B \in SO(2)$, $\xi$, $\eta \in P$,
where $P$ is the group of parallel translations on Euclidean plane.
We take exact representation of the group of proper movements
$x \mapsto Ax + \xi$,
$A \in SO(2), \xi \in R^2$,
in a form of $3 \times 3$ matrix
$$
\left(
\begin{array}{cc}
A & \xi \\
0 & 1
\end{array}
\right)
$$
proper orthonormal basis for Euclidean plane, which was earlier.

Lie algebra $L(E(2))$ of Lie group of proper movements
in Euclidean plane can be represented in the form
$L(E(2)) \cong so(2) {\oplus}_{M} R^2$, where the homomorphism
$M = d {\mu ( {\cdot}_{e} )}_{e} : so(2) \rightarrow Der R^2$;
$Der R^2$ is Lie algebra of the derivations of Lie algebra
$R^2$, and the action $\mu : SO(2) \rightarrow P$ is a rotation in
Euclidean plane.

Further we shall consider Lie algebra of Lie group
$E(2) \ast E(2)$.

We shall choose orthonormal basis in all
copies of $L(E(2))$ as follows:
$$
e_1 :=
\left( \left[ \begin {array}{cc} 0 & \frac{1}{\sqrt{2}} \\
\frac{-1}{\sqrt{2}} & 0 \end {array} \right], (0,0) \right),
$$
$$
e_2 :=
\left( \left[ \begin {array}{cc} 0 & 0 \\
0 & 0 \end {array} \right], (1,0) \right),
$$
$$
e_3 :=
\left( \left[
\begin {array}{cc} 0 & 0 \\
0 & 0 \end {array} \right], (0,1) \right).
$$

Since, as a vector space,
$L(E(2) \ast E(2)) \equiv L(E(2)) \oplus L(E(2))$,
then in Lie algebra $L(E(2) \ast E(2))$ naturally arises
orthonormal basis
$$
\{E_1,E_2,E_3,E_4,E_5,E_6\} =
\{e_1 \oplus 0, e_2 \oplus 0, e_3 \oplus 0,
0 \oplus e_1, 0 \oplus e_2, 0 \oplus e_3\}.
$$

On every copy of $L(E(2))$
Lie bracket of Lie algebra $L(E(2) \ast E(2))$
for the basis vectors
${\{E_i\}}_{i=1}^6$ accepts the values
(here we show only nonzero values):
$$
\begin{array}{lcr}
~[E_1,E_2] & = & -\frac{1}{\sqrt{2}} E_3, \vspace{4pt} \\
~[E_1,E_3] & = & \frac{1}{\sqrt{2}} E_2,
\end{array}
$$
$$
\begin{array}{lcr}
~[E_4,E_5] & = & -\frac{1}{\sqrt{2}} E_6, \vspace{4pt} \\
~[E_4,E_6] & = & \frac{1}{\sqrt{2}} E_5,
\end{array}
$$

Now we calculate the operators $L$ and $M$:
$$
L(A,\xi)(B,\eta) =
d {\lambda ( {(A,\xi)_{e} )}(B,\eta)}_{e} =
(0,A \eta),
$$
$$
M(B,\eta)(A,\xi) =
d {\mu ( {(B,\eta)_{e} )}(A,\xi)}_{e} =
(0,B \xi),
$$
where $A$, $B \in so(2)$, $\xi$, $\eta \in R^2$.

Thus, the operators $L$ and $M$
for the basis vectors
${\{E_i\}}_{i=1}^6$ accept the values
(here we show only nonzero values):
$$
\begin{array}{lcr}
L(E_4)(E_2) & = & -\frac{1}{\sqrt{2}} E_3, \vspace{4pt} \\
L(E_4)(E_3) & = & \frac{1}{\sqrt{2}} E_2,
\end{array}
$$
$$
\begin{array}{lcr}
M(E_1)(E_5) & = & -\frac{1}{\sqrt{2}} E_6, \vspace{4pt} \\
M(E_1)(E_6) & = & \frac{1}{\sqrt{2}} E_5,
\end{array}
$$

By the formulas (\ref{star}) of structural
constants of Lie algebra for twisted product of Lie groups
we receive the structural constants of Lie algebra $L(E(2) \ast E(2))$
(here we show only nonzero structural constants,
also granting the antisymmetry on the bottom indexes):
$$
\alpha_{12}^3=-\frac{1}{\sqrt{2}}, \alpha_{13}^2=\frac{1}{\sqrt{2}},
\alpha_{45}^6=-\frac{1}{\sqrt{2}}, \alpha_{46}^5=\frac{1}{\sqrt{2}};
$$
$$
\alpha_{15}^6=-\frac{1}{\sqrt{2}}, \alpha_{16}^5=\frac{1}{\sqrt{2}},
\alpha_{42}^3=-\frac{1}{\sqrt{2}}, \alpha_{43}^2=\frac{1}{\sqrt{2}};
$$

Then by the formula (\ref{sectcurv}) we receive a zero matrix
of sectional curvatures of the Lie group $E(2) \ast E(2)$.

From this we can easily find that the scalar curvature of the group
$E(2) \ast E(2)$ is equal $0$.
At the same time, as easily to check up,
all sectional curvatures of the initial group $E(2)$,
and consequently the scalar curvature, are equal to $0$.

This example is very interesting;
because it is known (see \cite[page 309]{M})
that the Euclidean group $E(2)$
is non-commutative, but admits a flat left invariant metric.
I.e. such metric that the sectional curvature $k_{ij}$
of Riemannian manifold $E(2)$ is identically zero.
\end{example}

\begin{example}{5}
As a new example we shall take, considered above in the example 4,
the twisted product of group $E(2)$ of proper
movements of the Euclidean plane by itself;
but we shall take another metric.

In generaly case, choosing some basis
$E_1$, ..., $E_n$ for the vector space $L(G)$,
it is easy to see that there is one and only one
Riemannian metric on $G$ so that vectors
$E_1$, ..., $E_n$ are orthonormal.

We shall choose basis in all copies of $L(E(2))$
as follows:
$$
e_1 := \left( \left[ \begin {array}{cc} 0 & \frac{1}{2} \\
\frac{-1}{2} & 0 \end {array} \right],
\left( \frac{1}{2},\frac{1}{2} \right) \right)
$$
$$
e_2 := \left( \left[ \begin {array}{cc} 0 & 0 \\
0 & 0 \end {array} \right], (1,0) \right)
$$
$$
e_3 := \left( \left[
\begin {array}{cc} 0 & 0 \\
0 & 0 \end {array} \right], (0,1) \right)
$$

In other words, for canonical metric the angle between $e_1$
and anyone of the vectors $e_2$, $e_3$
is $\pi / 3$.

Then there exist Riemannian metric on $E(2)$ so that vectors
$E_1$, ..., $E_3$ are orthonormal.
Further we shall consider this metric
on all copies of $E(2)$.

Since, as a vector space,
$L(E(2) \ast E(2)) \equiv L(E(2)) \oplus L(E(2))$,
then in Lie algebra $L(E(2) \ast E(2))$ naturally arises
orthonormal basis
$$
\{E_1,E_2,E_3,E_4,E_5,E_6\} =
\{e_1 \oplus 0, e_2 \oplus 0, e_3 \oplus 0,
0 \oplus e_1, 0 \oplus e_2, 0 \oplus e_3\}.
$$

On every copy of $L(E(2))$
Lie bracket of Lie algebra $L(E(2) \ast E(2))$
for the basis vectors
${\{E_i\}}_{i=1}^6$ accepts the values
(here we show only nonzero values):
$$
\begin{array}{lcr}
~[E_1,E_2] & = & -\frac{1}{2} E_3, \vspace{4pt} \\
~[E_1,E_3] & = & \frac{1}{2} E_2,
\end{array}
$$
$$
\begin{array}{lcr}
~[E_4,E_5] & = & -\frac{1}{2} E_6, \vspace{4pt} \\
~[E_4,E_6] & = & \frac{1}{2} E_5,
\end{array}
$$

The operators $L$ and $M$
for the basis vectors
${\{E_i\}}_{i=1}^6$ accept the values
(here we show only nonzero values):
$$
\begin{array}{lcr}
L(E_4)(E_1) & = & \frac{1}{4} E_2 - \frac{1}{4} E_3, \vspace{4pt} \\
L(E_4)(E_2) & = & -\frac{1}{2} E_3, \vspace{4pt} \\
L(E_4)(E_3) & = & \frac{1}{2} E_2,
\end{array}
$$
$$
\begin{array}{lcr}
M(E_1)(E_4) & = & \frac{1}{4} E_5 - \frac{1}{4} E_6, \vspace{4pt} \\
M(E_1)(E_5) & = & -\frac{1}{2} E_6, \vspace{4pt} \\
M(E_1)(E_6) & = & \frac{1}{2} E_5,
\end{array}
$$

By the formulas (\ref{star}) of structural
constants of Lie algebra for twisted product of Lie groups
we receive the structural constants of Lie algebra $L(E(2) \ast E(2))$
(here we show only nonzero structural constants,
also granting the antisymmetry on the bottom indexes):
$$
\alpha_{12}^3=-\frac{1}{2}, \alpha_{13}^2=\frac{1}{2},
\alpha_{45}^6=-\frac{1}{2}, \alpha_{46}^5=\frac{1}{2};
$$
$$
\alpha_{14}^2=-\frac{1}{4}, \alpha_{14}^3=\frac{1}{4},
\alpha_{14}^5=\frac{1}{4}, \alpha_{14}^6=-\frac{1}{4};
$$
$$
\alpha_{15}^6=-\frac{1}{2}, \alpha_{16}^5=\frac{1}{2},
\alpha_{42}^3=-\frac{1}{2}, \alpha_{43}^2=\frac{1}{2};
$$

Then by the formula (\ref{sectcurv}) we receive
a matrix of sectional curvatures
of the Lie group $E(2) \ast E(2)$:
$$
\left [\begin {array}{cccccc} 0&{\frac {1}{64}}&{\frac {1}{64
}}&\frac{-3}{16}&{\frac {1}{64}}&{\frac {1}{64}}\\\noalign{\medskip}{\frac
{1}{64}}&0&0&{\frac {1}{64}}&0&0\\\noalign{\medskip}{\frac {1}{64}}
&0&0&{\frac {1}{64}}&0&0\\\noalign{\medskip}\frac{-3}{16}&{\frac {1}{64}}&{
\frac {1}{64}}&0&{\frac {1}{64}}&{\frac {1}{64}}
\\\noalign{\medskip}{\frac {1}{64}}&0&0&{\frac {1}{64}}&0&0
\\\noalign{\medskip}{\frac {1}{64}}&0&0&{\frac {1}{64}}&0&0
\end {array}\right ]
$$

From this we can easily find that the scalar curvature of the group
$E(2) \ast E(2)$ is equal $-\frac{1}{8}$.
At the same time, as easily to check up,
all sectional curvatures of the initial group $E(2)$,
and consequently the scalar curvature, are equal to $0$.
\end{example}

\vspace{2\baselineskip}\noindent
{\bf Acknowledgements.}
The author is grateful to prof. V.N.Berestovskii for a kind attention
to the work.

\end{document}